\documentclass{article}
\usepackage{tabularx}
\usepackage{amsthm}
\newtheorem{corollary}{Corollary}

\usepackage{graphicx}
 \usepackage{amsmath} 
\newtheorem{proposition}{Proposition}

\usepackage[sglblindworkshop,final]{neurips_2025}
 \workshoptitle{MLxOR: Mathematical Foundations and Operational Integration of Machine Learning for Uncertainty-Aware Decision-Making}

\usepackage[utf8]{inputenc} 
\usepackage[T1]{fontenc}    
\usepackage{hyperref}       
\usepackage{url}            
\usepackage{booktabs}       
\usepackage{amsfonts}       
\usepackage{nicefrac}       
\usepackage{microtype}      
\usepackage{xcolor}         

\title{Overfitting in Adaptive Robust Optimization}

\author{
   Karl Zhu \\
   Operations Research Center \\
   Massachusetts Institute of Technology \\
     Cambridge, MA, USA \\
   \texttt{karlzhu@mit.edu} \\
  \And
  Dimitris Bertsimas \\
  Sloan School of Management \\
  Massachusetts Institute of Technology\\
  Cambridge, MA, USA \\
  \texttt{dbertsim@mit.edu} \\
}

\begin{document}

\maketitle

\begin{abstract}
Adaptive robust optimization (ARO) extends static robust optimization by allowing 
decisions to depend on the realized uncertainty ---
weakly dominating static solutions within the modeled uncertainty set.  However, ARO makes previous constraints that were independent of uncertainty now dependent, making it vulnerable to additional infeasibilities when realizations fall outside the uncertainty set. This phenomenon of adaptive policies being brittle is analogous to overfitting in machine learning. To mitigate against this, we propose assigning constraint-specific uncertainty set sizes, with harder constraints given stronger probabilistic guarantees. Interpreted through the overfitting lens, this acts as regularization: tighter guarantees shrink adaptive coefficients to ensure stability, while looser ones preserve useful flexibility. This view motivates a principled approach to designing uncertainty sets that balances robustness and adaptivity.
\end{abstract}

\section{Introduction}

Robust optimization (RO) is one of the main frameworks (along with stochastic programming) for decision-making under 
uncertainty. RO produces a \textit{robust} solution by requiring feasibility across all realizations within a specified uncertainty set $\mathcal U$. In this note, we consider a robust linear optimization problem with $m$ robust constraints and uncertainty on the right-hand side (RHS):
\begin{align*}
    \max_{\boldsymbol{x} \in \mathbb{R}^n} \quad & \boldsymbol{c}^\top \boldsymbol{x} \\
 \text{s.t.} \quad& \boldsymbol{a}_i^\top \boldsymbol{x} \le b_i + \boldsymbol{d}_i^\top \boldsymbol{u}, 
 \quad \forall \boldsymbol{u} \in \mathcal U, \quad \forall i \in [m], \\
 & \boldsymbol{x} \ge \boldsymbol{0},
\end{align*}
where for certain parameters we have $\boldsymbol{c}\in\mathbb{R}^n$,  and for each constraint $i$, we have  $\boldsymbol{a}_i\in\mathbb{R}^n, b_i\in\mathbb{R}, \boldsymbol{d}_i \in\mathbb{R}^p$. Finally, we have the perturbation parameter $\boldsymbol{u} \in\mathbb{R}^p$ which is uncertain but assumed to lie within $\mathcal{U}$.
Here, $\boldsymbol{x}$ is fixed before $\boldsymbol{u}$ is realized, yielding a \emph{static robust} solution. 

Adaptive robust optimization (ARO) extends this framework by allowing some 
decisions to adjust after observing $\boldsymbol{u}$ \cite{ben-tal_adjustable_2004}.
A common tractable restriction is to use \emph{affine decision rules}
\begin{equation}\label{eq:adr}
    \boldsymbol{x}(\boldsymbol{u}) = \boldsymbol{z} +  \boldsymbol{Vu}
\end{equation}
with decision variables $\boldsymbol{z}\in\mathbb{R}^n$ and $\boldsymbol{V}\in\mathbb{R}^{n\times p}$. 
Static solutions correspond to $\boldsymbol{V}=\boldsymbol{0}$, so adaptive policies always weakly dominate static ones within the set.

The size and geometry of $\mathcal U$ are chosen to balance feasibility and optimality. 
Common choices include the box ($\ell_\infty$ norm) \cite{soyster_technical_1973}, 
ellipsoid ($\ell_2$ norm) \cite{ben-tal_robust_1999}, and budget sets 
(intersection of $\ell_\infty$ and $\ell_1$ norms) \cite{bertsimas_price_2004}. 
A well-chosen uncertainty set generally should not contain all possible realizations of 
$\boldsymbol{u}$, as this would be overly conservative.
Therefore, when evaluating solutions in simulation, one should not draw samples directly 
from the uncertainty set, but rather from an underlying distribution. 

A key, under-discussed drawback of ARO is that constraints that were originally independent of 
uncertainty become dependent. While 
static robust solutions often remain feasible outside the modeled set, 
adaptive policies may fail catastrophically—for example, violating variable non-negativity constraints. In this sense, adaptive policies are 
\emph{brittle}: they achieve superior performance within $\mathcal U$ but 
generalize poorly outside it. This brittleness is directly analogous to 
\emph{overfitting} in machine learning, where models with higher flexibility 
fit training data well but perform poorly out-of-sample.

This perspective motivates a rethinking of how uncertainty sets are 
specified in ARO. Not all constraints are equally critical: hard 
constraints demand stronger guarantees, while softer constraints may 
tolerate limited violations. Interpreted through the overfitting lens, 
this leads to a \emph{regularization view}: tighter guarantees shrink 
adaptive coefficients to improve stability, while looser guarantees 
preserve flexibility. Although the equivalence of robustness and regularization has been studied extensively \cite{xu_robustness_nodate}, \cite{bertsimas_characterization_2018}, to our knowledge, there have been no similar discussions of this insight for ARO. In this note we illustrate brittleness using a 
renewable generation toy example, propose 
constraint-dependent uncertainty set sizes as a principled remedy for this brittleness, and show how robust counterparts (RC) impose regularization on adaptive coefficients.

\section{Brittleness of Adaptive Policies}

We illustrate brittleness with a toy model of production planning under 
renewables:
\begin{align*}
\max_{\boldsymbol{x},\boldsymbol{y},\boldsymbol{s}} \quad & x_1 + x_2 + y_1 + y_2 - 100s_1 - 100s_2 \\
\text{s.t.} \quad & x_1 + y_1 - s_1 \leq 2 + u_1, \quad \forall (u_1,u_2)\in\mathcal{U}, \\
                  & x_2 + y_2 - s_2 \leq 2 + u_2, \quad \forall (u_1,u_2)\in\mathcal{U}, \\
                  & x_1,x_2,y_1,y_2,s_1,s_2 \geq 0.
\end{align*}
Here a forecast of $2$ units of renewable power is available in each period. 
Decisions $x_1,y_1$ allocate this supply in period~1, and $x_2,y_2$ in 
period~2. Shortfalls are covered by grid imports $s_1,s_2$, which are costly 
and non-renewable. Uncertainty in renewable availability is represented by a 
budget uncertainty set
\[
\mathcal U = \{(u_1,u_2): \|\boldsymbol{u}\|_\infty \le \rho,\ \|\boldsymbol{u}\|_1 \le \Gamma\}.
\]
An optimal solution to the nominal problem (where $\mathcal{U} = \{\boldsymbol{0}\}$) is $x_i=y_i=1, \ s_i= 0$ for $i =1,2$ with an objective value of 4, but any shortfall 
in renewables forces costly grid usage.

An optimal solution to the static robust problem with $(\rho, \Gamma) = (1,1)$ is $x_i=1, y_i =s_i =0,$ for $ i=1,2$, with an objective value of 2. Notice that although the $l_1$ ball that defines $\mathcal{U}$ \textit{couples} the uncertainties across constraints sets in $\mathcal{U}$, the static robust problem in general cannot exploit this coupling \cite{marandi_when_2018}, giving an overly conservative solution that is equivalent to the uncertainty set being a box $
    \mathcal U = \{(u_1,u_2) : \|\boldsymbol{u}\|_\infty \le \rho\}.$

In contrast, ARO is able to exploit this coupling \cite{bertsimas_benefit_2024}. Suppose we introduce a restricted affine decision rules that made the flows adapt to the current period of uncertainty:
\[
x_i(u_i) = z^x_i + V^x_{i,i} u_i, \quad 
y_i(u_i) = z^y_i + V^y_{i,i} u_i, \quad i =1,2.
\] 
An optimal adaptive policy with $(\rho, \Gamma)=(1,1)$ is $x_i(u_i)=1,\ y_i(u_i)=1+u_i$, for $i=1,2$ with an objective value of 3 -- significantly outperforming the static 
policy within $\mathcal U$. However, the adaptive policy is brittle: the 
non-negativity constraints, originally independent of $u$, now depend on $u$ 
and may be violated. For instance, if $u_1<-1$, then $y_1(u_1)<0$, which is 
physically infeasible. Shortages can be met with imports, but negative 
flows cannot be implemented. 

A natural remedy is to assign larger uncertainty radii to hard constraints 
(e.g., non-negativity) and smaller ones to softer constraints (e.g., grid imports). 
For example, with $(\rho,\Gamma)=(2,2)$ on the non-negativity constraints only, an optimal adaptive 
solution is $x_i=y_i=1+\tfrac12 u_i$, for $i=1,2$ and maintains an objective value of 3 and remains feasible for 
all $u\ge -2$. In our example, only a non-negative RHS of renewable supply is physically sensible, so $u_1,u_2$ naturally both have a lower bound of $-2$. So the enlarged set recovers the guaranteed feasibility previously enjoyed by the  static policy while retaining the benefits from ARO and uncertainty coupling.

\section{Probabilistic and Deterministic Guarantees}
From the toy example we saw that a uniform uncertainty set $\mathcal U$ is not 
appropriate in ARO, since many constraints that were previously independent of 
uncertainty become dependent once adaptivity is introduced. This motivates the 
use of \emph{constraint-specific uncertainty sets} $\mathcal U_i$ for each 
constraint $i \in [m]$, sized according to its criticality. We distinguish two 
types:
\begin{itemize}
    \item \textbf{Hard constraints} (e.g. flow non-negativity): 
    must be satisfied in all realizations. These require deterministic 
    guarantees, with uncertainty sets that fully cover the support of 
    $\boldsymbol{u}$.
    \item \textbf{Softer constraints} (e.g., renewable allocations with 
    non-renewable grid imports): can tolerate limited violations if backup resources 
    exist. For these, probabilistic guarantees are appropriate, where 
    ellipsoidal or budget sets provide explicit guarantees while preserving ARO's flexibility.
\end{itemize}

We now provide quantitative prescriptions. First, we recall probabilistic 
guarantees under Gaussian and distribution-free assumptions from the RO literature \cite{bertsimas_probabilistic_2021}, suited for softer 
constraints. Then, we present deterministic RCs for bounded and 
semi-bounded supports, which are appropriate for hard constraints that must hold 
almost surely.

\paragraph{Probabilistic guarantees.}
For each constraint $i\in[m]$, we may require
\begin{equation}\label{eq:prob-guarantee}
\mathbb P\!\left[\boldsymbol{a}_i^\top \boldsymbol{x}(\boldsymbol{u})
\;\le\; b_i + \boldsymbol{d}_i^\top \boldsymbol{u}\right]\ \ge\ 1-\varepsilon.    
\end{equation}
Different uncertainty sets $\mathcal U_i$ provide different sufficient conditions 
for \eqref{eq:prob-guarantee}.

\begin{itemize}
\item \textbf{Gaussian case.}  
If $\boldsymbol{u}\sim\mathcal N(\boldsymbol{\mu},\Sigma)$ and 
$\boldsymbol{x}(\boldsymbol{u})$ is affinely adaptive (\ref{eq:adr}), then 
\eqref{eq:prob-guarantee} holds iff
\[
\boldsymbol{a}_i^\top \boldsymbol{z} 
+ (\boldsymbol{V}^\top \boldsymbol{a}_i - \boldsymbol{d}_i)^\top \boldsymbol{\mu} 
+ \rho_{1-\varepsilon}\,\bigl\|\Sigma^{1/2}(\boldsymbol{V}^\top \boldsymbol{a}_i - \boldsymbol{d}_i)\bigr\|_2
\;\le\; b_i,
\]
where $\rho_{1-\varepsilon}$ is the $(1-\varepsilon)$-quantile of the standard normal. 
Equivalently, this is the RC with ellipsoidal uncertainty set
\[
\mathcal U_i=\left\{\boldsymbol{u}:\ 
\bigl\|\Sigma^{-1/2}(\boldsymbol{u}-\boldsymbol{\mu})\bigr\|_2 \le \rho_{1-\varepsilon}\right\}.
\]

\item \textbf{Distribution-free case.}  
If $u_1,\dots,u_p$ are independent, zero-mean, with support $[-1,1]$, then:
\begin{itemize}
\item \emph{Ball-box uncertainty set.}  
For $\mathcal U_i=\{\boldsymbol{u}:\|\boldsymbol{u}\|_\infty\le 1,\ \|\boldsymbol{u}\|_2\le \rho\}$, the RC is
\begin{equation}\label{eq:rc-ball}
\boldsymbol{a}_i^\top \boldsymbol{z} + \rho\,\|\boldsymbol{V}^\top \boldsymbol{a}_i - \boldsymbol{d}_i - \boldsymbol{y}_i\|_2 + \|\boldsymbol{y}_i\|_1 \;\le\; b_i,
\quad \boldsymbol{y}_i\in\mathbb{R}^p,
\end{equation}
which ensures
$\mathbb P[\boldsymbol{a}_i^\top \boldsymbol{x}(\boldsymbol{u})\le b_i + \boldsymbol{d}_i^\top \boldsymbol{u}]
\ge 1-\exp(-\rho^2/2)$. Choosing $\rho=\sqrt{2\ln(1/\varepsilon)}$ yields 
\eqref{eq:prob-guarantee}.
\item \emph{Budget uncertainty set.}  
For $\mathcal U_i =\{\boldsymbol{u}: \|\boldsymbol{u}\|_\infty\le 1,\ \|\boldsymbol{u}\|_1\le \Gamma \}$, the RC is
\[
\boldsymbol{a}_i^\top \boldsymbol{z} + \Gamma\,\|\boldsymbol{V}^\top \boldsymbol{a}_i - \boldsymbol{d}_i - \boldsymbol{y}_i\|_\infty + \|\boldsymbol{y}_i\|_1 \;\le\; b_i,
\quad \boldsymbol{y}_i\in\mathbb{R}^p,
\]
which ensures $\mathbb P[\boldsymbol{a}_i^\top \boldsymbol{x}(\boldsymbol{u})\le b_i + \boldsymbol{d}_i^\top \boldsymbol{u}]
\ge 1-\exp(-\Gamma^2/(2p))$. Choosing $\Gamma=\sqrt{2\ln(1/\varepsilon)}\sqrt{p}$ yields 
\eqref{eq:prob-guarantee}.
\end{itemize}
\end{itemize}

\paragraph{Deterministic guarantees.}
For hard constraints that must hold without exception,
probabilistic guarantees are insufficient. In these cases, one must construct an uncertainty 
set that fully contains the support of $\boldsymbol{u}$ and enforce feasibility deterministically. 
This leads to robust counterparts over polyhedral box uncertainty sets, whose derivation follows 
standard duality arguments (see \cite{bertsimas_robust_2006} for a tutorial). We state the 
results directly below.

\begin{proposition}[Bounded support]\label{prop:bounded}
Suppose $u_k$ has support $[L_k,U_k] \ \forall k \in [p]$, and $\boldsymbol{x}(\boldsymbol{u})$ is affinely adaptive (\ref{eq:adr}). Then the following RC satisfies $\boldsymbol{a}_i^\top \boldsymbol{x}(\boldsymbol{u})\le b_i + \boldsymbol{d}_i^\top \boldsymbol{u}$ w.p. 1:
\[
\boldsymbol{a}_i^\top \boldsymbol{z} + \boldsymbol{U}^\top \boldsymbol{\beta}  -  \boldsymbol{L}^\top \boldsymbol{\alpha} \ \le b_i, \quad  \boldsymbol{\beta} - \boldsymbol{\alpha} = \boldsymbol{V}^\top \boldsymbol{a}_i - \boldsymbol{d}_i, \quad \boldsymbol{\alpha,\beta} \ge \boldsymbol{0},\]
where $\boldsymbol{L} = (L_k)_{k=1}^p, \boldsymbol{U} = (U_k)_{k=1}^p$ and $\boldsymbol{\alpha}, \boldsymbol{\beta} \in \mathbb{R}^p$.
\end{proposition}

\begin{corollary}[Semi-bounded support]
Under the same setup as Proposition~\ref{prop:bounded}, but instead if $u_k$ is unbounded from below (resp. above), for each $k \in [p]$, the RC from Proposition~\ref{prop:bounded} along with $\alpha_k = 0 \ (\text{resp. } \beta_k =0 )$ satisfies 
$\boldsymbol{a}_i^\top \boldsymbol{x}(\boldsymbol{u})\le b_i + \boldsymbol{d}_i^\top \boldsymbol{u}$ w.p. 1.
\end{corollary}
\section{A Regularization Perspective}

\begin{figure}[h]
    \centering
    \includegraphics[width=0.65\linewidth]{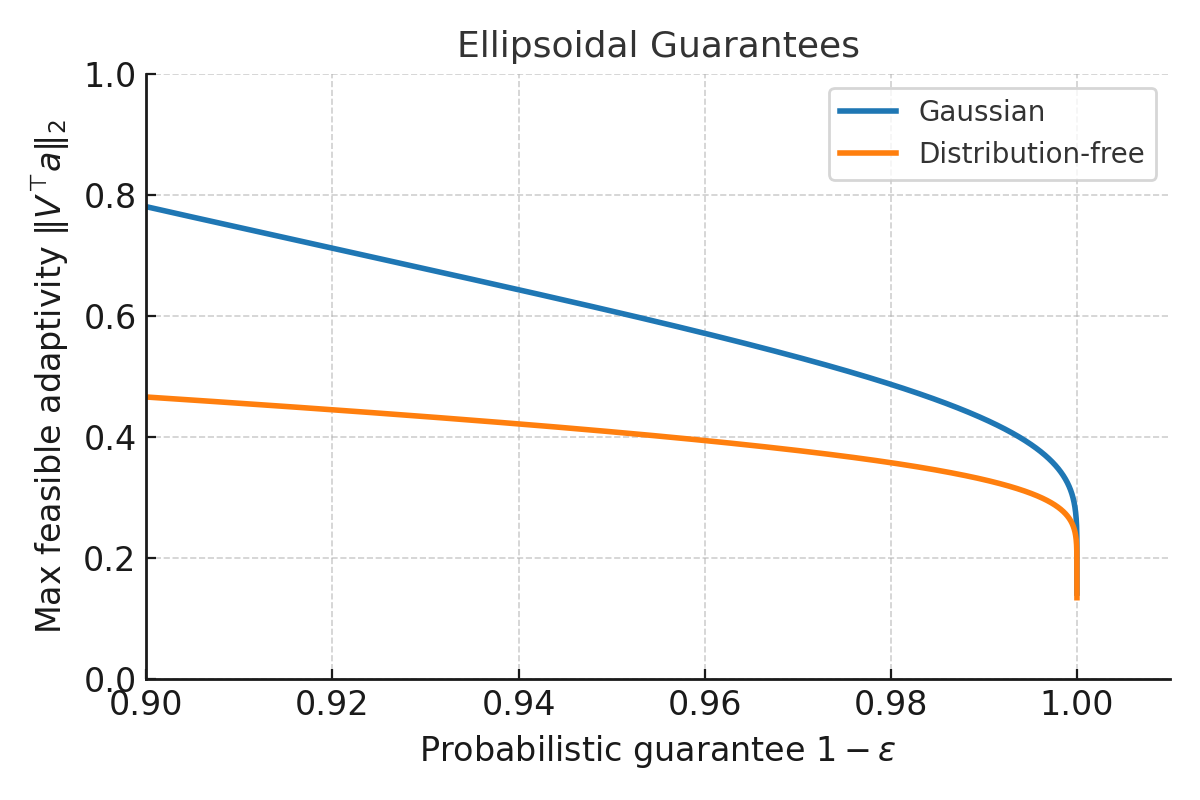}
\caption{Ellipsoidal guarantees: maximum feasible adaptivity 
$\|\boldsymbol{V}^\top \boldsymbol{a}\|_2$ as a function of the probabilistic guarantee $1-\varepsilon$. 
Under Gaussian assumptions, the underlying probabilistic bound is tighter, yielding less conservative 
feasibility regions. The distribution-free guarantee is looser, leading to more conservative regularization. 
In both cases, higher guarantees correspond to stronger regularization, shrinking adaptive flexibility 
and preventing brittle solutions.}
    \label{fig:regularization-path}
\end{figure}

We now interpret robust constraints through the lens of \emph{regularization}, 
providing a mitigation for the brittleness of ARO in line with the 
overfitting analogy. Notice the RCs can be 
rearranged into explicit norm bounds on the adaptive coefficients. 
For example, with $\boldsymbol{y}_i = \boldsymbol{0}$, (\ref{eq:rc-ball}) can be rewritten as
\[
\|\boldsymbol{V}^\top \boldsymbol{a}_i - \boldsymbol{d}_i\|_2 \;\le\; \frac{\boldsymbol{b}_i - \boldsymbol{a}_i^\top \boldsymbol{z}}{\rho}, 
\qquad \rho > 0.
\]
For illustration, we set $\boldsymbol{d}_i=\boldsymbol{0}$ and normalize $\boldsymbol{b}_i - \boldsymbol{a}_i^\top \boldsymbol{z} = 1$, 
and plot the resulting upper bound on $\|\boldsymbol{V}^\top \boldsymbol{a}_i\|_2$ against the 
probabilistic guarantee $1-\varepsilon$. 
Figure~\ref{fig:regularization-path} compares the bounds derived under 
a Gaussian assumption with those obtained from the 
distribution-free guarantee.  

In this view, the probabilistic guarantee $1-\varepsilon$ serves as a \emph{regularization parameter}: requiring stronger guarantees (smaller $\varepsilon$) tightens the constraint on $\boldsymbol{V}$, limiting adaptivity. In the limit as $\varepsilon \to 0$, we recover the static policy with $\boldsymbol{V} = \boldsymbol{0}$, eliminating brittleness. The resulting $\ell_2$-norm bound directly parallels ridge regression in machine learning, where stronger regularization shrinks coefficients to improve stability. This mirrors the {bias--variance tradeoff}: greater adaptivity yields better in-set performance but poorer generalization, while stronger regularization reduces flexibility but stabilizes out-of-set behavior.

\section{Discussion and Conclusion}

Our analysis suggests several lessons for modeling with ARO. 
\begin{itemize}
    \item \textbf{Simulate beyond the set.} Evaluate performance on 
    distributions whose support exceeds the modeled uncertainty set, 
    just as machine learning models must be tested out-of-sample.
    \item \textbf{Differentiate constraints.} Treat some as \emph{soft}, 
    where limited violations are tolerable (e.g., supply--demand balance), 
    and others as \emph{hard}, where violations are catastrophic 
    (e.g., flow non-negativity).
    \item \textbf{Use constraint-dependent uncertainty set sizes.} Assign larger sets 
    to hard constraints, effectively regularizing their adaptive 
    coefficients more strongly, while smaller sizes suffice for softer ones.
\end{itemize}
Overall, ARO provides flexibility but can overfit to the specified
uncertainty sets, producing brittle out-of-set behavior. Interpreting robust 
counterparts as implicit regularization motivates constraint-dependent uncertainty sets 
as a principled way to balance adaptivity and stability.
\clearpage

\bibliography{aro_overfit}

\bibliographystyle{unsrt}

\end{document}